\documentclass[12pt,twoside]{article}
\usepackage{amsmath, amssymb, amsthm}
\newtheorem{theorem}{Theorem}

\newtheorem{lemma}{Lemma}

\newtheorem{definition}{Definition}

\font\fourteenb=cmb10 at 14pt
\begin{document}
\vspace*{-1.0cm}\noindent \copyright
 Journal of Technical University at Plovdiv\\[-0.0mm]\
\ Fundamental Sciences and Applications, Vol. 5, 1997\\[-0.0mm]
\textit{Series A-Pure and Applied Mathematics}\\[-0.0mm]
\ Bulgaria, ISSN 1310-8271\\[+1.2cm]
\font\fourteenb=cmb10 at 14pt
\begin{center}

   {\bf \LARGE Free interpolation of families of Cauchy-Stiltjes integrals and their multipliers
   \\ \ \\ \large Peyo Stoilov}
\end{center}

\

\

\begin{abstract}
The present paper contains a generalization of some interpolation
theorems of S. A. Vinogradov [3].
\end{abstract}

\section{Introduction}
Let $\Bbb{D}\,$ denotes the unit disk in the complex plane and $\Bbb{T}\,$-
the unit circle. Let $M\,$ be the Banach space of all complex-valued Borel
measures on $\Bbb{T}\,$ with the usual variation norm. For $\alpha>0$, let
${\cal F}_{\alpha}$ denote the family of functions $g\,$ for which there exists
$\mu\in M\,$ such that
\begin{equation}
g(z)=\int\limits_{T}\frac{1}{(1-\overline{\xi}z)^{\alpha}}\,d\mu(\xi),
\quad z\in\Bbb{D}.
\end{equation}
We note that ${\cal F}_{\alpha}\,$ is a Banach space with the natural norm
\[
\Vert g\Vert_{{\cal F}_{\alpha}}=\inf\left\{\Vert \mu\Vert :\mu\in M
\mbox{ such that (1) holds}\;\right\}.
\]
\begin{definition}
Suppose that $f\,$ is holomorphic in $\Bbb{D}$. Then $f\,$ is called a
{\bf multiplier} of ${\cal F}_{\alpha}\,$ if $g\in{\cal F}_{\alpha}
\Rightarrow fg\in{\cal F}_{\alpha}$.
\end{definition}
Let $m_{\alpha}\,$ denote the set of multipliers of ${\cal F}_{\alpha}\,$
and
\[
\Vert f\Vert_{m_{\alpha}}=\sup\left\{\Vert fg\Vert_{{\cal F}_{\alpha}}
:\Vert g\Vert_{{\cal F}_{\alpha}}\leq 1\right\}.
\]
We shall note some results we shall use further.
\def\thetheorem{A}
\begin{theorem}
Let $\alpha>0$. Then $f\in m_{\alpha}\,$ if and only if
\[
f(z)\frac{1}{(1-\overline{\xi}z)^{\alpha}}\in{\cal F}_{\alpha}
\]
for $|\xi|=1\,$ and there
is a positive constant $A\,$ such that
\[
\Vert f(z)
\frac{1}{(1-\overline{\xi}z)^{\alpha}}\Vert_{{\cal F}_{\alpha}}\leq A
\]
for $|\xi|=1$.
\end{theorem}
\def\thetheorem{B}
\begin{theorem}
Let $0<\alpha<\beta$. Then $\;m_{\alpha}\subset m_{\beta}\subset H^{\infty}\,$ and
\[
\Vert f\Vert_{\infty}\leq
\Vert f\Vert_{m_{\beta}}\leq \Vert f\Vert_{m_{\alpha}}.
\]
\end{theorem}
\def\thetheorem{C}
\begin{theorem}
If the function $f\,$ is analytic in $\Bbb{D}\,$ and
\[
\Omega(f)\stackrel{def}{=}\int\limits_0^1\,\int\limits_{-\pi}^{\pi}
|f(re^{\imath\, \theta})|(1-r)^{\alpha-1}\,d\theta\,dr<\infty
\]
for some $\alpha>0$, then $f\in {\cal F}_{\alpha+1}\,$ and
$\Vert f\Vert_{{\cal F}_{\alpha+1}}\leq \frac{\alpha}{2\pi}\,\Omega(f)$.
\end{theorem}
\def\thetheorem{D}
\begin{theorem}
Let $\alpha>0$. Then $f\in{\cal F}_{\alpha}\,$ if and only if
$f'\in{\cal F}_{\alpha+1}$. In addition
\begin{equation}
\frac{1}{\alpha}\Vert f'\Vert_{{\cal F}_{\alpha+1}}\leq
\Vert f\Vert_{{\cal F}_{\alpha}}\leq |f(0)|+\frac{2}{\alpha}
\Vert f'\Vert_{{\cal F}_{\alpha+1}}.
\end{equation}
\end{theorem}
Theorems A and B were proved in [2]. From the proof of Theorem A it follows
that
\[
\Vert f\Vert_{m_{\alpha}}=\sup\left\{\left\Vert \frac{f(z)}
{(1-\overline{\xi}z)^{\alpha}}\right\Vert_{{\cal F}_{\alpha}}:|\xi|=1\right\}.
\]
Theorem C is Lemma 1 in [4]. Theorem D was proved in [1]. On the inequalities (2), see [4, Lemma 2].

In this paper some interpolation theorems due to S. A. Vinogradov for $m_1\,$ and ${\cal F}_1\:$
[2,Theorems 9 and 12] are generalized for $m_{\alpha}\,$ and $
{\cal F}_{\alpha},\;\alpha>0$.
\section{Free interpolation in $m_{\alpha}$}
\begin{definition}
We say that a sequence $a=\{a_k\}_{k \geq 1}\subset {\Bbb D}\,$ satisfies
Newman--Carleson condition (condition ($N-C$)) if
\begin{equation}
\delta(a)=\inf\left\{\prod\limits_{k \neq n}\left\vert\frac{a_k-a_n}
{1-\overline{a}_k a_n}\right\vert :\;n=1,2,\dots\right\}>0.
\end{equation}
\end{definition}
\begin{lemma}
Let a sequence $a=\{a_k\}_{k \geq 1}\subset {\Bbb D}\,$ satisfy condition
($N-C$) and condition
\begin{equation}
\sigma_{\alpha}(a)=\sup\limits_{|\xi|=1}\sum\limits_{k\geq 1}\left(\frac
{1-|a_k|^2}{|1-\overline{\xi}a_k|}\right)^{\alpha}<\infty,
\end{equation}
for some $0<\alpha \leq 1$. Then for each sequence $x=\{x_k\}_{k\geq 1}\in
\ell^{\infty}$, there is a function $f\,$ in $m_{\alpha}$, such that
\begin{description}
\item[(a)]$f(a_k)=x_k\;,\quad k=1,2,\dots$ ;
\item[(b)]$\Vert f\Vert_{m_{\alpha}}\leq c\left(\frac{\sigma_{\alpha}(a)}
{\delta(a)}\right)^2 \Vert x\Vert_{\ell^{\infty}}$,
\end{description}
where $c>0\,$ is a constant depending only on $\alpha$.
\end{lemma}
{\it Proof}.\hspace{.3cm}The case $\alpha=1\,$ was proved in [3]. Let $0<\alpha
<1\,$ and $x=\{x_k\}_{k\geq 1}\in\ell^{\infty}$. Consider the function that
is used in [3]:
\[
f(z)=\sum\limits_{n\geq 1}\frac{1-|a_n|^2}{1-\overline{a}_n z}\,\frac{B_n(z)}
{B_n(a_n)}\,x_n,
\]
where
\[
B_n(z)=\prod\limits_{
\begin{array}{c}
k=1\\
k\neq n
\end{array}
}^{\infty}\frac{z-a_k}{1-\overline{a}_k z}\,\frac{|a_k|}{a_k}\,,\quad
n=1,2,\dots.
\]
Since the sequence $a\,$ satisfies the condition ($N-C$), then from Lemma 9
([3]) it follows that
\[
\frac{B_n(z)}{1-\overline{a}_n z}=\frac{a_n}{|a_n|}\sum\limits_{k\geq 1}
\frac{|a_k|}{a_k}\,\frac{1-|a_k|^2}{1-\overline{a}_k a_n}\,\frac{1}
{\overline{B_k(a_k)}}\,\frac{1}{1-\overline{a}_k z}.
\]
This implies that
\begin{eqnarray*}
f(z)=\sum\limits_{n\geq 1}\frac{1-|a_n|^2}{B_n(a_n)}\,\frac{a_n}{|a_n|}x_n
\sum\limits_{k\geq 1}
\frac{|a_k|}{a_k}\,\frac{1-|a_k|^2}{1-\overline{a}_k a_n}\,\frac{1}
{\overline{B_k(a_k)}}\,\frac{1}{1-\overline{a}_k z}=\\[10pt]
\sum\limits_{k\geq 1}\frac{1}{1-\overline{a}_k z}\,
\sum\limits_{n\geq 1}\frac{(1-|a_n|^2)(1-|a_k|^2)}{B_n(a_n)
\overline{B_k(a_k)}}\,\frac{a_n|a_k|}{|a_n|a_k}\,\frac{x_n}
{1-\overline{a}_k a_n}=\\[10pt]
\sum\limits_{k\geq 1}\frac{y(a_k)}{1-\overline{a}_k z},
\end{eqnarray*}
where
\[
y(a_k)\stackrel{def}{=}\sum\limits_{n\geq 1}\frac{(1-|a_n|^2)(1-|a_k|^2)}{B_n(a_n)
\overline{B_k(a_k)}}\,\frac{a_n|a_k|}{|a_n|a_k}\,\frac{x_n}
{1-\overline{a}_k a_n}.
\]
Making use of (3) and (4) we estimate $y(a_k)$:
\begin{equation}
|y(a_k)|\leq\Vert x\Vert_{\ell^{\infty}}\,\frac{1-|a_k|^2}{\delta^2}
\sum\limits_{n\geq 1}\frac{1-|a_n|^2}{|1-\overline{a}_k a_n|}\leq
\frac{\sigma_1}{\delta^2}\Vert x\Vert_{\ell^{\infty}}\,(1-|a_k|^2),
\end{equation}
where $\delta=\delta(a)\,,\; \sigma_1=\sigma_1(a)$. We note that $\sigma_1\leq \sigma_{\alpha}\,$ for $0<\alpha<1$.

Further, Theorem A will be used to prove that $f\in m_{\alpha}$, or equivalently
\begin{equation}
\Vert f\Vert_{m_{\alpha}}=\sup\{\Vert g\Vert_{{\cal F}_{\alpha}}:
|\xi|=1\}<\infty,
\end{equation}
where
\[
g(z)\stackrel{def}{=}f(z)\left(\frac{1}{1-\overline{\xi}z}\right)^{\alpha}
\,,\;|\xi|=1.
\]
Theorem D implies that it suffices to show that
\[
g'(z)=f'(z)\frac{1}{(1-\overline{\xi}z)^{\alpha}}+\alpha\overline{\xi}f(z)
\frac{1}{(1-\overline{\xi}z)^{\alpha+1}}
\]
belongs to ${\cal F}_{\alpha+1}\,$ and there is a constant $B\,$ such that
\begin{equation}
\Vert g'\Vert_{{\cal F}_{\alpha+1}}\leq B\; \mbox{for all} \;\xi\in {\Bbb
T}.
\end{equation}
Since $0<\alpha< 1\,$ and $\sigma_1(a)\leq \sigma_{\alpha}(a)<\infty$, then
by the theorem of Vinogradov $f\in m_1\,$ and
\[
\Vert f\Vert_{m_1}\leq 2\left(\frac{\sigma_1}{\delta}\right)^2
\Vert x\Vert_{\ell^{\infty}}.
\]
Thus theorem B implies that $f\in m_{\alpha+1}\,$ and consequently
\[
\Vert\alpha\overline{\xi}f(z)\frac{1}{(1-\overline{\xi}z)^{\alpha+1}}
\Vert_{{\cal F}_{\alpha+1}}\leq \alpha \Vert f\Vert_{m_{\alpha+1}}\leq
\alpha \Vert f\Vert_{m_1}\leq 2\alpha\left(\frac{\sigma_1}{\delta}\right)^2
\Vert x\Vert_{\ell^{\infty}}.
\]
For the proof of (7) it remains to be shown that
\[
f'(z)\frac{1}{(1-\overline{\xi}z)^{\alpha}}\stackrel{def}{=}h(z)\in
{\cal F}_{\alpha+1}
\]
too. Since
\[
f'(z)=\sum\limits_{k\geq1}\frac{y(a_k)}{(1-\overline{a}_k z)^2}\,
\overline{a}_k,
\]
from (5) it follows that
\[
|f'(z)|\leq\frac{\sigma_1}{\delta^2}
\Vert x\Vert_{\ell^{\infty}}\sum\limits_{k\geq1}\frac{1-|a_k|^2}{|1-
\overline{a}_k z|^2}.
\]
Then
\begin{eqnarray*}
\Omega(h)=\int\limits_0^1\,\int\limits_{-\pi}^{\pi}|h(re^{\imath\theta})|
(1-r)^{\alpha-1}\,d\theta\,dr\leq\\
\frac{\sigma_1}{\delta^2}\Vert x\Vert_{\ell^{\infty}}
\int\limits_0^1\,\int\limits_{-\pi}^{\pi}\sum\limits_{k\geq1}\frac
{1-|a_k|^2}{|1-\overline{a}_k re^{\imath\theta}|^2}\,\frac
{(1-r)^{\alpha-1}}{|1-\overline{\xi} re^{\imath\theta}|^{\alpha}}\,d\theta
\,dr=\\
\frac{\sigma_1}{\delta^2}\Vert x\Vert_{\ell^{\infty}}\sum\limits_{k\geq1}\,
\int\limits_0^1\,\int\limits_{-\pi}^{\pi}\frac
{1-|a_k|^2}{|1-\overline{a}_k re^{\imath\theta}|^2}\,\frac
{(1-r)^{\alpha-1}}{|1-\overline{\xi} re^{\imath\theta}|^{\alpha}}\,d\theta
\,dr.
\end{eqnarray*}
In [4] (Lemma 13) the inequality
\[
\int\limits_0^1\,\int\limits_{-\pi}^{\pi}\frac{(1-r)^{\alpha-1}}
{|1-\overline{z} re^{\imath\theta}|^2
|1-\overline{\xi} re^{\imath\theta}|^{\alpha}}\,d\theta\,dr\leq E\frac
{(1-|z|)^{\alpha-1}}{|1-\overline{\xi}z|^{\alpha}}
\]
was proved for $|\xi|=1\,,\;|z|<1\,$ and a constant $E>0\,$ depending only
on $\alpha$. This inequality implies that
\[
\Omega(h)\leq E \frac{\sigma_1}{\delta^2}\Vert x\Vert_{\ell^{\infty}}\sum\limits_{k\geq1}\frac{(1-|a_k|^2)(1-|a_k|)^
{\alpha-1}}{|1-\overline{\xi}a_k|^{\alpha}}\leq
2E\left(\frac{\sigma_{\alpha}}{\delta}\right)^2\Vert x\Vert_{\ell^{\infty}}
\]
and by Theorem C $\,h(z)\in {\cal F}_{\alpha+1}$, and
\[
\Vert h\Vert_{{\cal F}_{\alpha+1}}\leq \alpha E
\left(\frac{\sigma_{\alpha}}{\delta}\right)^2\Vert x\Vert_{\ell^{\infty}}.
\]
This completes the proof of the fact that $g'\in{\cal F}_{\alpha+1}\,$ and
\[
\Vert g'\Vert_{{\cal F}_{\alpha+1}}\leq \alpha (2+E)
\left(\frac{\sigma_{\alpha}}{\delta}\right)^2\Vert x\Vert_{\ell^{\infty}}.
\]
Since by Theorem D
\[
\Vert g\Vert_{{\cal F}_{\alpha}}\leq |g(0)|+\frac{2}{\alpha}
\Vert g'\Vert_{{\cal F}_{\alpha+1}}
\]
and
\[
|g(o)|=|f(0)|\leq \Vert f\Vert_{\infty}\leq \Vert f\Vert_{m_1}\leq
2\left(\frac{\sigma_{\alpha}}{\delta}\right)^2\Vert x\Vert_{\ell^{\infty}},
\]
it follows from (6) that
\[
\Vert f\Vert_{m_{\alpha}}\leq C \left(\frac{\sigma_{\alpha}}{\delta}\right)^2\Vert x\Vert_{\ell^{\infty}},
\]
where $C=6+2E$.

This proves (b). The equalities (a)
\[
f(a_k)=x_k\,,\;k=1,2,\dots\,,
\]
are obvious. $\;\;\Box$
\section{Free interpolation in ${\cal F}_{\alpha}$}
\begin{lemma}
Let a sequence $a=\{a_k\}_{k \geq 1}\subset {\Bbb D}\,$ satisfy condition
($N-C$) and condition
\begin{equation}
\sigma_{\alpha}(a)=\sup\limits_{|\xi|=1}\sum\limits_{k\geq 1}\left(\frac
{1-|a_k|^2}{|1-\overline{\xi}a_k|}\right)^{\alpha}<\infty,
\end{equation}
for some $0<\alpha \leq 1$. Then
\begin{description}
\item[(a)] If $g\in {\cal F}_{\alpha}$, then $(g(a_k)(1-|a_k|^2)^{\alpha}
)_{k\geq 1}\,\in \ell^1\,$;
\item[(b)] For each sequence $x=\{x_k\}_{k\geq 1}\in
\ell^{1}$, there is a function $g\,$ in ${\cal F}_{\alpha}$, such that
\[
g(a_k)(1-|a_k|^2)^{\alpha}=x_k\quad,\,k=1,2,\dots.
\]
\end{description}
\end{lemma}
Proof.\hspace{.3cm}The case $\alpha =1\,$ was proved in [3]. Let $0<\alpha<
1$.\\
Proof of {\bf (a)}. Let $g\in {\cal F}_{\alpha}\,$ and
\[
g(z)=\int\limits_{T}\frac{1}{(1-\overline{\xi}z)^{\alpha}}\,d\mu(\xi).
\]
Then
\[
\sum\limits_{k\geq 1}|g(a_k)|(1-|a_k|^2)^{\alpha}=
\sum\limits_{k\geq 1}\left|\int\limits_{T}\left(\frac
{1-|a_k|^2}{|1-\overline{\xi}a_k|}\right)^{\alpha}\,d\mu(\xi)\right|
\leq \\ \Vert\mu\Vert\ \sigma_{\alpha}(a)<\infty.
\]
Proof of {\bf (b)}. Let $x=\{x_k\}_{k\geq 1}\in\ell^{1}\,$ and
\[
g(z)=\sum\limits_{n\geq 1}\frac{B_n(z)}{(1-\overline{a}_n z)^{\alpha}}\,\frac{x_n}
{B_n(a_n)}\,,
\]
where
\[
B_n(z)=\prod\limits_{
\begin{array}{c}
k=1\\
k\neq n
\end{array}
}^{\infty}\frac{z-a_k}{1-\overline{a}_k z}\,\frac{|a_k|}{a_k}\,,\quad
n=1,2,\dots.
\]
Since $|B_n(a_k)|\geq \delta(a)>0$, to prove that $g\in {\cal F}_{\alpha}\,$
it suffices to show that
\[
g_n(z)=\frac{B_n(z)}{(1-\overline{a}_n z)^{\alpha}}\in {\cal F}_{\alpha}\;\;
\mbox{and}\; \sup\limits_n \Vert g_n \Vert_{{\cal F}_{\alpha}}<\infty.
\]
From the proof of Theorem 5 [4] it follows that
\[
\Vert B_n \Vert_{m_{\alpha}}\leq D\,,
\]
where the constant $D\,$ depends only $\alpha\,$ and $\sigma_{\alpha}(a)$.
Since $\Vert \frac{1}{(1-zw)^{\alpha}} \Vert_{{\cal F}_{\alpha}}=1\,$ for
$|w|\leq1\,$(Lemma 6, [4]), then
\[
\Vert g_n \Vert_{{\cal F}_{\alpha}}\leq \Vert B_n \Vert_{m_{\alpha}}\leq D.
\]
Therefore $g\in {\cal F}_{\alpha}$. The equalities
\[
g(a_k)(1-|a_k|^2)^{\alpha}=x_k\quad,\,k=1,2,\dots
\]
are obvious. $\;\;\Box$

\section{Remarks}
Condition (4), i.e. $\sigma_{\alpha}(a)<
\infty\,$ was introduced in [4]. For $\alpha=1\,$ this condition is the
uniform Frostman condition.

We also note that in [5] some interpolation theorems for $m_{\alpha}\,$ and ${\cal F}_{\alpha}\,$ were proved in the case when a sequence
$a=\{a_k\}_{k\geq1}\,$ satisfies condition ($N$-$C$) and all points of the
sequence lie in a single Stolz region.
\vspace{2cm}\\

\centerline{{\bf REFERENCES}}

\begin{enumerate}
\item {\sc T. H. Mac Gregor}, Analytic and univalent functions
with integral representations involving complex measures, {\it
Indiana Univ. Math. J.} {\bf 36} (1987), 109-130. \item {\sc R. A.
Hibschweiler and T. H. Mac Gregor}, Multipliers of families of
Cauchy-Stieltjes transforms, {\it Trans. Amer. Math. Soc.} {\bf
331} (1992), 377-394. \item {\sc S. A. Vinogradov}, Properties of
multipliers of Cauchy-Stieltjes integrals and some factorization
problems for analytic functions, {\it Amer. Math. Soc. Transl.} 2,
{\bf 115} (1980), 1-32. \item {\sc D. J. Hallenbeck, T. H. Mac
Gregor and K. Samotij}, Fractional Cauchy transforms, inner
functions and multipliers, {\it Proc. London Math. Soc.} (3) {\bf
72} (1996), 157-187. \item {\sc P. Stoilov}, On some interpolation
theorems for the multipliers of the Cauchy-Stieltjes type
integrals, {\it J. Tech. Univ. at Plovdiv, Fundam.  Sci. Appl. ,}
{\bf 2} (1996), 33-40.

\

\

\noindent
{\small Department of Mathematics\\
        Technical University\\
        25, Tsanko Dijstabanov,\\
        Plovdiv, Bulgaria\\
        e-mail: peyyyo@mail.bg}
        \end{enumerate}
\end{document}